\documentclass[12pt]{article}
\usepackage{amssymb} 

\def\d{\delta}

\def\g{\gamma}

\def\i{\iota}

\def\o{\omega}

\def\u{\upsilon}

\def\O{\Omega}

\chardef\tempcat=\the\catcode`\@
\catcode`\@=11
\def\cyracc{\def\u##1{\if \i##1\accent"24 i
    \else \accent"24 ##1\fi }}
\newfam\cyrfam

\DeclareFontFamily{OT1}{msb}{}{}
\DeclareFontShape{OT1}{msb}{m}{n}
 {  <5> <6> <7> <8> <9> <10> gen * msbm
      <10.95><12><14.4><17.28><20.74><24.88>msbm10}{}
\DeclareMathAlphabet{\bubble}{OT1}{msb}{m}{n}

\def\bR{{\mathbb R}}

\def\bC{{\mathbb C}}

\newfont{\goth}{eufm10 scaled \magstep1}

\newfont{\mcal}{eusm10 scaled \magstep1}

\def\square{\kern1pt\vbox
            {\hrule height 0.6pt\hbox{\vrule width 0.6pt\hskip 3pt
 \vbox{\vskip 6pt}\hskip 3pt\vrule width 0.6pt}\hrule height 0.6pt}\kern1pt}

\def\ra{\rightarrow}

\def\wt{\widetilde}

\newtheorem{Th}{Theorem}
\newtheorem{Prop}{Proposition} 
\newtheorem{Cor}{Corollary} 
\newtheorem{Lem}{Lemma}
\newtheorem{Def}{Definition} 

\def\bt{\begin{Th}}
\def\et{\end{Th}}
\def\bp{\begin{Prop}}
\def\ep{\end{Prop}}
\def\bc{\begin{Cor}}
\def\ec{\end{Cor}}
\def\bl{\begin{Lem}}
\def\el{\end{Lem}}
\def\bd{\begin{Def}}
\def\ed{\end{Def}}

\def\pf{\noindent{\it Proof:\ }}
\def\qed{\hfill\square}
\def\n{\nabla}

\def\be{\begin{equation}}
\def\ee{\end{equation}}

\def\arr{\begin{array}{rlll}}
\def\ea{\end{array}}
\def\bea{\begin{eqnarray}}
\def\eea{\end{eqnarray}}  
\def\bean{\begin{eqnarray*}}
\def\eean{\end{eqnarray*}}  

\catcode`@=11
\@addtoreset{equation}{section}
\catcode`@=12

\begin{document}
\title{ A holomorphic representation formula\\
       for parabolic hyperspheres} 

\author{Vicente Cort\' es \thanks{e-mail: vicente@math.uni-bonn.de \newline
This work was supported by SFB 256 ``Nichtlineare partielle 
Differentialgleichungen'' (Universit\"at Bonn).} \\
Mathematisches Institut \\
Universit\"at Bonn\\
Beringstra\ss e 1 \\
D-53115 Bonn}
\date{\today}
\maketitle
\begin{abstract} \noindent  
A holomorphic representation formula
       for special parabolic hyperspheres is given. 


  \end{abstract}

\noindent
\section*{Introduction}
It was noticed by Blaschke that parabolic spheres in affine 
$3$-space admit parametrisations in terms of 
holomorphic functions (of one variable). This is related to the fact  
that the Monge-Amp\`ere equation governing parabolic hyperspheres 
is completely 
integrable in dimension two, a fact already known to Monge. 
The purpose of this note is to derive an explicit formula describing 
special parabolic hyperspheres in affine $(2n+1)$-space in terms of 
a holomorphic function of $n$ variables. 
\section{Special parabolic affine hyperspheres}
Let me briefly recall the notion of {\it special} parabolic affine hypersphere 
\cite{BC}. We consider  $\bR^{m+1}$ as affine space with standard 
connection denoted by $\widetilde{\n}$ and parallel volume 
form ${\rm vol}$. A hypersurface is given by an immersion 
$\varphi : M \rightarrow \bR^{m+1}$ of an $m$-dimensional connected 
manifold. We assume that $M$ admits a transversal vector field $\xi$ and 
that $m>1$. 
This induces
on $M$ the volume form $\nu = {\rm vol} 
(\xi , \ldots )$, a torsionfree connection
$\n$, a quadratic covariant tensor field $g$, an endomorphism field $S$ 
(shape tensor) and a 
one-form $\theta$ such that
\bea \widetilde{\n}_XY &=& \n_XY + g(X,Y)\xi\, ,\nonumber\\ 
 \widetilde{\n}_X\xi &=& SX + \theta (X)\xi \, .\label{GWEqu}\eea
Let us call the data $(\n ,g,S,\theta )$ the 
{\bf Gau{\ss}-Weingarten data} induced by the transversal vector field $\xi$. 
We will assume that $g$ is nondegenerate and, hence, is a 
pseudo-Riemannian metric on $M$. This condition does not depend
on the choice of $\xi$.  
According to Blaschke \cite{Blaschke}, 
once the orientation of $M$ is fixed, there is a unique
choice of transversal vector field $\xi$ such that
$\nu$ coincides with the metric volume form ${\rm vol}^g$ and $\n \nu = 0$. 
This particular choice of  transversal vector field is called the 
{\bf affine normal} and the corresponding geometric data $(g,\n )$ are called
{\bf Blaschke metric} and {\bf induced connection}. 
Notice that, for the affine normal, $\theta = 0$ and 
$S$ is computable from $(g,\n )$ (Gau{\ss} equations).  
Henceforth we use always the affine normal as transversal
vector field.  
\bd \label{specialDef} 
The hypersurface $\varphi : M \rightarrow \bR^{m+1}$ is called
a {\bf parabolic} (or improper) {\bf hypersphere} if the affine
normal is parallel, $\widetilde{\n}\xi = 0$. 
It is called {\bf special} if there exists an  
almost complex structure $J$ on $M$ which is skew symmetric with respect
to the Blaschke metric $g$ and such that the 2-form 
$\omega := g(J \cdot , \cdot )$ is parallel with respect to the induced 
connection $\n$. Such an almost complex structure $J$ is called 
{\bf compatible}.  
\ed 
Notice that $\widetilde{\n}\xi = 0 \Leftrightarrow S = 0 
\Leftrightarrow \n$ is flat. It was proven in \cite{BC} that any  
parabolic two-dimensional sphere with (positive or negative) 
definite Blaschke metric is special and that a compatible
almost complex structure on a special parabolic hypersphere is 
necessarily integrable. In fact, we proved the following stronger result, 
Theorem \ref{BCThm} below. 

Recall that a  {\bf special K\"ahler manifold}  $(M,J,g, \n )$
is a (pseudo-)K\"ahler manifold $(M,J,g)$ endowed with a 
flat torsionfree connection $\n$ such that $\n J$ is symmetric and 
$\n \o = 0$, where $\omega = g(J \cdot , \cdot )$ is the  K\"ahler 
form. 
\bt \cite{BC} \label{BCThm} 
Let $\varphi : M \rightarrow \bR^{m+1}$ be a special   
parabolic hypersphere with Blaschke metric $g$, induced connection $\n$, 
compatible almost complex structure $J$ and canonical
two-form $\omega = g(J \cdot , \cdot )$. Then $(M,J,g, \n )$ is a
special K\"ahler manifold. 
Conversely, any simply connected special K\"ahler manifold
$(M,J,g, \n )$ admits an immersion $\varphi : M \rightarrow \bR^{m+1}$,
which is a special parabolic hypersphere with Blaschke metric $g$, 
induced connection $\n$ and compatible almost complex structure $J$.
The immersion $\varphi$ is unique up to a unimodular affine transformation
of $\bR^{m+1}$. 
\et  
\section{The holomorphic representation formula} 
It was proven in \cite{ACD} that any special K\"ahler manifold is locally
defined by a holomorphic function, as follows. Let $F$ be a holomorphic
function on a domain (i.e.\ a connected open set) 
$U \subset \bC^n$ such that the (real) matrix 
\be \label{nondegEqu} {\rm Im}\, \partial^2 F \quad \mbox{is invertible,}\ee
where $\partial^2F$ denotes the 
holomorphic Hessian of $F$. Let us denote by $M_F \subset T^*\bC^n$ 
the image of the holomorphic section $dF : U \ra T^*U \subset T^*\bC^n$.
It is a complex Lagrangian submanifold with respect to the standard 
complex symplectic structure $\O = \sum dz^i\wedge dw_i$, where
$(z^1, \ldots ,z^n, w_1, \ldots , w_n)$ are canonical coordinates
of  $T^*\bC^n$. We denote its complex structure by $J$.  
Using the nondegeneracy condition (\ref{nondegEqu}), it is shown in 
\cite{ACD} that the Hermitian form $\gamma  := \sqrt{-1}
\O (\cdot , \bar{\cdot})$ 
is nondegenerate on $M_F$ and, hence, induces a (pseudo-) K\"ahler metric 
$g = {\rm Re}\, \g |_{M_F}$. It is also shown that  a 
torsionfree connection $\n$ on 
$M_F$ can be defined by the condition that the real parts 
$x^i :=  {\rm Re}\, z^i$ and $y_j := {\rm Re}\, w_j$ are $\n$-affine functions 
on $M_F$. In fact, it is shown that 
$(x^1, \ldots ,x^n, y_1, \ldots , y_n)$ is a (real) local
coordinate system near any point of $M_F$ and that the K\"ahler form
$\o = g (\cdot , J\cdot )$ is expressed by the formula $\o = 
2\sum dx^i\wedge dy_i$ on $M_F$. 

\bt \cite{ACD} \label{ACDThm} Let $F$ be a holomorphic function 
satisfying the nondegeneracy condition (\ref{nondegEqu}) 
on a domain $U \subset \bC^n$. Then $(M_F,J,g,\n )$, defined above,  
is a special K\"ahler manifold and any 
special K\"ahler manifold is locally of this form. 
\et 
It is noticed in \cite{BC} that combining Theorem \ref{BCThm} and 
Theorem \ref{ACDThm} we can associate a parabolic hypersphere
to any holomophic function $F$ defined on a simply connected domain 
$U \subset \bC^n$ and satisfying the nondegeneracy condition
(\ref{nondegEqu}). However, the proof of Theorem \ref{BCThm} makes use of 
the Fundamental Theorem of affine differential geometry \cite{DNV} 
(the generalisation
of Radon's theorem \cite{R} to higher dimensions) 
and does not involve any explicit
parametrisation of the immersion   
$\varphi : M \rightarrow \bR^{2n+1}$ realising 
a simply connected special K\"ahler
manifold $(M,J,g,\n )$ of real dimension $2n$ as parabolic hypersphere.  
The aim is now to provide an explicit formula, in 
terms of the holomorphic function $F$, for the realisation
of $(M_F,J,g,\n )$ as a parabolic hypersphere 
$\varphi_F : M_F \cong U \rightarrow \bR^{2n+1}$. 

We will not restrict ourselves to functions $F$ defined on 
simply connected domains $U \subset \bC^n$. More generally, we consider a 
`multivalued' function defined on an arbitrary domain $U \subset 
\bC^n$. Or, in other words, a (univalued) function defined on some 
Riemann domain $\wt{U}$ over $U$. A {\bf Riemann domain} over $U$ is a 
holomorphic (unramified) covering $\pi : \wt{U} \ra U$.  Any holomorphic
function $F$ on $\wt{U}$ defines a holomorphic Lagrangian immersion 
\be \label{phiEqu} 
\phi : \wt{U} \ra  T^*U \subset T^*\bC^n\, , \quad \phi (p) := 
dF \circ (\pi_{\ast}|T_p\wt{U})^{-1} \, ,\quad p \in \wt{U}
\, .\ee
Let us denote by $J$ the complex structure of $\wt{U}$. 
Pulling back the canonical coordinates of $T^*\bC^n$ to $\wt{U}$ 
we obtain holomorphic functions
\[ \tilde{z}^i := \phi^* z^i \quad \mbox{and} \quad \tilde{w}_j := \phi^* w_j 
\] 
on $\wt{U}$. The holomorphic functions $\tilde{z}^i$ form a local 
holomorphic coordinate system near any point of $\wt{U}$. We use the 
compact notation 
\[ \tilde{z} := (\tilde{z}^1, \ldots ,\tilde{z}^n)\, , 
\quad F_{\tilde{z}} = (F_{\tilde{z}^1}, \ldots ,
F_{\tilde{z}^n}) = (\frac{\partial F}{\partial \tilde{z}^1}, \ldots ,
\frac{\partial F}{\partial \tilde{z}^n})\, ,\quad   F_{\tilde{z}}\tilde{z} =
\sum F_{\tilde{z}^k}\tilde{z}^k\quad \mbox{etc.} 
\] 
Let $\partial^2F$
be the Hessian of $F$ with respect to (the flat torsionfree 
holomorphic connection defined by) the coordinate system $\tilde{z}$. 
We call
$F$ {\bf nondegenerate} if ${\rm Im}\, \partial^2 F$ is invertible. 
Then, as before,  
$g := {\rm Re}\, \phi^*\g$ is a pseudo-K\"ahler metric and we can define
a flat torsionfree connection $\n$ by the condition that the functions
$\tilde{x}^i :=  {\rm Re}\, \tilde{z}^i$ and 
$\tilde{y}_j := {\rm Re}\, \tilde{w}_j$ are\linebreak[4] 
$\n$-affine functions 
on $\wt{U}$. We also put $\tilde{u}^i := {\rm Im }\, \tilde{z}^i$ and 
$\tilde{v}_j :=  {\rm Im }\, \tilde{w}_j$. 
Let us abbreviate $M(F) := (\wt{U},J,g,\n )$ and define an immersion
$\varphi_F : \wt{U} \ra \bR^{2n+1}$ by the formula
\bea \varphi_F  &:=& ({\rm Re}\, \tilde{z}, {\rm Re}\, F_{\tilde{z}}, 
2{\rm Im}\, F - 2({\rm Re}\, F_{\tilde{z}}){\rm Im}\, \tilde{z}) 
\label{varphiEqu} \\
&=& (\tilde{x}^1,\ldots, \tilde{x}^n, \tilde{y}_1, 
\ldots , \tilde{y}_n, 
2{\rm Im}\, F - 2\sum({\tilde{y}_k})\tilde{u}^k)\nonumber\, .
\eea
\pagebreak[3]
\begin{samepage}   
\bt \label{HRThm} Let $F$ be a nondegenerate 
holomorphic function defined on a Riemann domain 
$\wt{U}$. Then $M(F) = (\wt{U},J,g,\n )$, defined above, 
is a special K\"ahler manifold with K\"ahler form $\o = 
g(\cdot ,J \cdot ) = 2\sum d\tilde{x}^i\wedge d\tilde{y}_i$. The immersion  
$\varphi_F : \wt{U} \ra \bR^{2n+1}$ defined by 
(\ref{varphiEqu}) is, with respect to the volume form 
${\rm vol} := 2^n\det$ on 
$\bR^{2n+1}$, a special parabolic hypersphere with 
affine normal $\xi = \partial_{2n+1}$, Blaschke metric $g$, 
induced connection $\n$ and compatible 
almost complex structure $J$. It is unique up to unimodular 
affine transformations of $\bR^{2n+1}$. 
\et    
\end{samepage} 
\pf The first statement is a slight generalisation of the first part of 
Theorem~\ref{ACDThm}, with essentially the same proof. The uniqueness of
$\varphi_F$ follows, as in the proof of Theorem~\ref{BCThm},  
from the uniqueness statement of the Fundamental Theorem of affine differential
geometry. It suffices to prove that $\varphi_F$ is a parabolic 
hypersphere with Blaschke metric $g$, induced connection $\n$ and compatible 
almost complex structure $J$. Let us compute the Gau{\ss}-Weingarten
data $(\n^v,g^v,S^v,\theta^v)$, see (\ref{GWEqu}),  induced by the 
transversal vector field $v = \partial_{2n+1}$ (the `vertical' vector field). 
It is immediate that $S^v = 0$ and $\theta^v = 0$.  
We compute $\n^v$ and $g^v$ for the 
coordinate vector fields
\[ \partial_{\tilde{x}^i} = 
\partial_i + \frac{\partial f}{\partial \tilde{x}^i}\partial_{2n+1}\, ,\quad 
 \partial_{\tilde{y}_j}  = 
\partial_{n+j} + \frac{\partial f}{\partial \tilde{y}_j}\partial_{2n+1}
\, ,
\] 
where 
\[ f :=  2{\rm Im}\, F - 2({\rm Re}\, F_{\tilde{z}}){\rm Im}\, \tilde{z} =
 2{\rm Im}\, F -2\sum \tilde{y}_k\tilde{u}^k 
\]   
is the 
last component of $\varphi_F$. 
The covariant derivatives with respect to the connection $\wt{\n}$ of 
$\bR^{2n+1}$ are given by  
\[ \wt{\n}_{\partial_{\tilde{x}^i}} \partial_{\tilde{x}^j} = 
\frac{\partial^2 f}{\partial_{\tilde{x}^i}\partial_{\tilde{x}^j}}v\, ,\quad 
\wt{\n}_{\partial_{\tilde{x}^i}} \partial_{\tilde{y}_j} = 
\wt{\n}_{\partial_{\tilde{y}_j}}\partial_{\tilde{x}^i} = 
\frac{\partial^2 f}{\partial_{\tilde{x}^i}\partial_{\tilde{y}_j}}v\, ,\quad \wt{\n}_{\partial_{\tilde{y}_i}} 
\partial_{\tilde{y}_j} = 
\frac{\partial^2 f}{\partial_{\tilde{y}_i}\partial_{\tilde{y}_j}}v\, .\] 
This shows that the coordinate vector fields  $\partial_{\tilde{x}^i}$ 
and $\partial_{\tilde{y}_j}$ are parallel for the connection $\n^v$, so
it coincides with $\n$. Now $\theta = 0$ implies $\n \nu^v = \n^v \nu^v = 
0$ for the volume form $\nu^v = {\rm vol} (v , \ldots )$. 
Moreover, we see that $g^v = {\rm Hess}^\n(f) = \n^2 f$.\\ 
{\bf Claim 1} $g^v = g$.\\ 
The claim, to be proven below, 
implies that $v$ is the affine normal and, hence, that $g^v=g$ is 
the Blaschke metric. Let us see why. The Riemannian volume of the 
(pseudo-)K\"ahler manifold $M(F)$ with K\"ahler form 
$\o = g(J \cdot , \cdot ) = 2\sum_{i=1}^n d\tilde{x}^i\wedge d\tilde{y}_i$ 
is given by 
\[ {\rm vol}^g = 
(-1)^{n(n-1)/2}\frac{\o^n}{n!} = 2^n 
d\tilde{x}^1\wedge \ldots \wedge d\tilde{x}^n\wedge  
d\tilde{y}_1\wedge \ldots \wedge d\tilde{y}_n = 2^n \det (v, \ldots ) = 
\nu^v\, ,\] 
if we choose the orientation defined by $\nu^v$. (Notice that $d\tilde{x}^1\wedge \ldots \wedge d\tilde{x}^n\wedge  
d\tilde{y}_1\wedge \ldots \wedge d\tilde{y}_n = (-1)^{n(n-1)/2}
d\tilde{x}^1\wedge d\tilde{y}_1 \wedge \ldots \wedge d\tilde{x}^n\wedge 
d\tilde{y}_n$.) This shows that 
${\rm vol}^g$ coincides with the $\n$-parallel volume form $\nu^v$. 
So $v$ is the affine normal and, hence,     
$\varphi_F$ is a parabolic hypersphere with 
Blaschke metric $g$ and induced connection $\n$. {}The fact that
$M(F)$ is a special K\"ahler manifold entails that $J$ is 
skew symmetric with respect to $g$ and $\n \o = 0$. Therefore,
the complex structure
$J$ is compatible, in the sense of Definition \ref{specialDef},  
with the data $(g,\n )$. 

\noindent 
It remains to prove Claim 1. For the calculations we will use the next lemma.
\bl \label{OLemma} The partial derivatives of the functions $\tilde{u}^i$ and 
$\tilde{v}_j$ on $M(F)$ with respect to the $\nabla$-affine coordinates 
$(\tilde{x}^1,\ldots, \tilde{x}^n, \tilde{y}_1, 
\ldots , \tilde{y}_n)$ satisfy the following equations:
\bean  
&& \sum_k(\tilde{u}^k_{\tilde{x}^i}(\tilde{v}_k)_{\tilde{y}_j} - 
\tilde{u}^k_{\tilde{y}_j}(\tilde{v}_k)_{\tilde{x}^i}) = \d^j_i \, 
  \, ,
\\
&& \sum_k \tilde{u}^k_{\tilde{x}^i}(\tilde{v}_k)_{\tilde{x}^j} = \sum_k 
\tilde{u}^k_{\tilde{x}^j}(\tilde{v}_k)_{\tilde{x}^i} \, 
 ,\\
&&  \sum_k \tilde{u}^k_{\tilde{y}_i}(\tilde{v}_k)_{\tilde{y}_j} = \sum_k 
\tilde{u}^k_{\tilde{y}_j}(\tilde{v}_k)_{\tilde{y}_i}\,\\ 
&&   \tilde{u}^i_{\tilde{x}^j} = - (\tilde{v}_j)_{\tilde{y}_i}\, , 
\\
&&   \tilde{u}^i_{\tilde{y}_j} =  \tilde{u}^j_{\tilde{y}_i} \, , 
\\
&&   (\tilde{v}_i)_{\tilde{x}^j} =  (\tilde{v}_j)_{\tilde{x}^i} \,  .  
\eean  
\el 

\pf Pulling back the symplectic form $\O$ of $T^*\bC^n$ by means of 
the Lagrangian immersion $\phi : M(F) \ra T^*\bC^n$ defined in (\ref{phiEqu}), 
we obtain the equation $\phi^* \O = 0$. Decomposing it into real
and imaginary parts yields the lemma. \qed

\noindent 
Let us return to the proof of Theorem \ref{HRThm}. First we observe that
\[  \partial_{\tilde{x}^i} {\rm Im}\, F =   {\rm Im}\, \partial_{\tilde{x}^i} 
F = \sum_j {\rm Im}\, (\frac{\partial \tilde{z}^j}{\partial \tilde{x}^i}
\partial_{\tilde{z}^j} F) = \sum_j {\rm Im}\, ((\d_i^j + 
\sqrt{-1} \tilde{u}^j_{\tilde{x}^i})F_{\tilde{z}^j}) =  
\tilde{v}_i +  \sum_j \tilde{u}^j_{\tilde{x}^i}\tilde{y}_j\] 
and 
\[ \partial_{\tilde{y}_j} {\rm Im}\, F = {\rm Im}\, \partial_{\tilde{y}_j} F 
=  \sum_k {\rm Im}\, (\frac{\partial \tilde{z}^k}{\partial \tilde{y}_j}
\partial_{\tilde{z}^k} F) = \sum_k {\rm Im}\, (\sqrt{-1}
\tilde{u}^k_{\tilde{y}_j} F_{\tilde{z}_k}) = \sum_k 
\tilde{u}^k_{\tilde{y}_j}\tilde{y}_k\, . 
\]  
The second derivatives of ${\rm Im}\, F$ are now easily 
computed with the help of
Lemma \ref{OLemma}:
\[  \partial^2_{\tilde{x}^i\tilde{x}^j} {\rm Im}\, F = 
(\tilde{v}_i)_{\tilde{x}^j} + \sum_k \tilde{u}^k_{\tilde{x}^i\tilde{x}^j}
\tilde{y}_k\, ,\quad \partial^2_{\tilde{x}^i\tilde{y}_j} 
{\rm Im}\, F = \sum_k  \tilde{u}^k_{\tilde{y}_j\tilde{x}^i}\tilde{y}_k 
\, ,\quad 
 \partial^2_{\tilde{y}_i\tilde{y}_j} {\rm Im}\, F = \tilde{u}^i_{\tilde{y}_j}
+ \sum_k \tilde{u}^k_{\tilde{y}_i\tilde{y}_j}\tilde{y}_k\, .
\] 
Using this and Lemma \ref{OLemma} one can now evaluate $g^v = \n^2 f$: 
\bea g^v(\partial_{\tilde{x}^i},\partial_{\tilde{x}^j}) &=& 
\partial^2_{\tilde{x}^i\tilde{x}^j} f = 2 ((\tilde{v}_i)_{\tilde{x}^j}
+ \sum_k \tilde{u}^k_{\tilde{x}^i\tilde{x}^j}\tilde{y}_k) - 2 
\sum_k \tilde{u}^k_{\tilde{x}^i\tilde{x}^j}\tilde{y}_k = 
2 (\tilde{v}_i)_{\tilde{x}^j}\, ,\nonumber \\
g^v(\partial_{\tilde{x}^i},\partial_{\tilde{y}_j}) &=& 
\partial^2_{\tilde{x}^i\tilde{y}_j} f = 2 
\sum_k \tilde{u}^k_{\tilde{y}_j\tilde{x}^i} \tilde{y}_k 
- 2 (\tilde{u}^j_{\tilde{x}^i} + \sum_k \tilde{u}^k_{\tilde{y}_j
\tilde{x}^i}\tilde{y}_k) = - 2 \tilde{u}^j_{\tilde{x}^i}\, ,\label{symmEqu}\\
g^v(\partial_{\tilde{y}_i},\partial_{\tilde{y}_j}) &=& 
\partial^2_{\tilde{y}_i\tilde{y}_j} f =  2(\tilde{u}^i_{\tilde{y}_j}
+ \sum_k \tilde{u}^k_{\tilde{y}_i\tilde{y}_j}\tilde{y}_k) - 
2(\tilde{u}^j_{\tilde{y}_i} + \tilde{u}^i_{\tilde{y}_j} +
\sum_k \tilde{u}^k_{\tilde{y}_i\tilde{y}_j}\tilde{y}_k)  
= -2\tilde{u}^i_{\tilde{y}_j}\, .\nonumber
\eea 
Notice that in virtue of (\ref{symmEqu}) we have:
\be (\tilde{u}^i)_{\tilde{x}^j} = (\tilde{u}^j)_{\tilde{x}^i}\, . 
\label{symmEqu'} 
\ee 
Let us compare this with $g$. The simplest way to compute $g$ is using
the fact that $g = \o \circ J$, where $\o = 2\sum dx^i\wedge dy_i$ is 
the K\"ahler form and we consider $g$ and $\o$ as isomorphisms $TM \ra T^*M$
(insertion of a vector in the first argument). 
It is easier to work with the inverse metric $g^{-1} = J^{-1} \circ
\o^{-1} = - J \circ \o^{-1} = \o^{-1} \circ J^*$. Notice that
\[ \o^{-1} = \frac{1}{2} \sum \partial_{\tilde{y}_i}\wedge
\partial_{\tilde{x}^i}\, ,\quad 
J^*d\tilde{x}^i = -d\tilde{u}^i\quad \mbox{and} \quad 
J^*d\tilde{y}_j = -d\tilde{v}_j\, .\]
Let us evaluate $g^{-1}$ with the help of these formulas and Lemma 
\ref{OLemma}:
\bean g^{-1} (d\tilde{x}^i,d\tilde{x}^j)  
&=& -\o^{-1}(d\tilde{u}^i,d\tilde{x}^j) 
=  -\frac{1}{2}\tilde{u}^i_{\tilde{y}_j}\, ,\\
g^{-1} (d\tilde{x}^i,d\tilde{y}_j) &=& -\o^{-1}(d\tilde{u}^i,d\tilde{y}_j) = 
\frac{1}{2}\tilde{u}^i_{\tilde{x}^j} \, ,\\
g^{-1} (d\tilde{y}_i,d\tilde{y}_j) &=& -\o^{-1}(d\tilde{v}_i,d\tilde{y}_j)
= \frac{1}{2}(\tilde{v}_i)_{\tilde{x}^j} \,.
\eean
Comparing with the formulas for $g^v$ and using 
Lemma \ref{OLemma} and (\ref{symmEqu'}) this proves that
$g^{-1}g^v = {\rm id}$ and, hence, that $g = g^v$. This completes
the proof of Claim 1 and Theorem~\ref{HRThm}.  
\qed 

\end{document}